\documentclass[english]{article}
\usepackage{epstopdf}
\usepackage[T1]{fontenc}
\usepackage[latin9]{inputenc}
\usepackage{geometry}
\usepackage{float}
\usepackage{mathtools}
\usepackage{amsmath}
\usepackage{amsthm}
\usepackage{amssymb}
\usepackage{graphicx}
\usepackage{subfig}
\usepackage{setspace}
\usepackage{xcolor}
\usepackage{romanbar}
\usepackage{cite}

\makeatletter

\providecommand{\tabularnewline}{\\}

\numberwithin{equation}{section}
\numberwithin{figure}{section}
\theoremstyle{plain}

\theoremstyle{remark}

\makeatother

\makeatletter
\renewcommand*\env@matrix[1][\arraystretch]{%
	\edef\arraystretch{#1}%
	\hskip -\arraycolsep
	\let\@ifnextchar\new@ifnextchar
	\array{*\c@MaxMatrixCols c}}

\makeatother

\usepackage{babel}
\providecommand{\remarkname}{Remark}
\providecommand{\theoremname}{Theorem}

\newcommand{\tL}{\theta^L_{ij}}
\newcommand{\tR}{\theta^R_{ij}}

\newcommand{\normal}{\mathbf{n}}
\newcommand{\uvec}{u}
\newcommand{\vel}{V}
\newcommand{\xy}[2]{\mathbf{X}_{#1,#2}}
\newcommand{\pt}{\mathbf{X}}
\newcommand{\ptij}{\mathbf{X}_{ij}}

\begin{document}

\title{A second-order accurate semi-Lagrangian method for convection-diffusion equations with interfacial jumps}
\author{Hyuntae Cho, Yesom Park, and Myungjoo Kang}
\maketitle
\begin{abstract}
	In this paper, we present a second-order accurate finite-difference method for solving convection-diffusion equations with interfacial jumps on a moving interface. The proposed method is constructed under a semi-Lagrangian framework for convection-diffusion equations; a novel interpolation scheme is developed in the presence of jump conditions. Combined with a second-order ghost fluid method \cite{cho2019second}, a sharp capturing method with a first-order local truncation error near the interface and second-order truncation error away from the interface is developed for the convection-diffusion equation. In addition, a level-set advection algorithm is presented when the velocity gradient jumps across the interface. Numerical experiments support the conclusion that the proposed methods for convection-diffusion equations and level-set advection are necessary for the second-order convergence solution and the interface position.
\end{abstract}

\section{Introduction}

In this article, we consider the following convection-diffusion equation:
\begin{equation}\label{eq:PDE}
\begin{aligned}
\rho \left(\uvec_t+\vel \cdot\nabla \uvec\right)-\nabla\cdot\left(\mu \nabla \uvec\right) & =f\ \text{\ensuremath{\text{in}}\ \ensuremath{\Omega\textbackslash\ensuremath{\Gamma_t}\times (0,T)}},\\
\uvec & =g\ \ensuremath{\text{on}}\ensuremath{\ \partial\Omega\times(0,T)}, \\
\uvec & =h\ \ensuremath{\text{on}}\ensuremath{\ \Omega\times\{0\}},
\end{aligned}
\end{equation}
where $\Omega$ is a bounded domain in $\mathbb{R}^{n}$ $\left(n=1,2\right)$ with boundary $\partial\Omega$ and $\Gamma_t$ is a codimension-1 moving interface
that divides $\Omega$ into disjoint subdomains $\Omega^{+}_t$ and
$\Omega^{-}_t$ at time $t\in \left(0,T\right)$. Furthermore, the solution for \eqref{eq:PDE} satisfies the jump {conditions}
\begin{equation}\label{eq:jmp}
\begin{aligned}
[\mu \frac{\partial \uvec}{\partial \normal}]&=b   \text{\ \ensuremath{\text{on}}\ensuremath{\ \Gamma_t}}, \\
[\uvec]&=a   \text{\ \ensuremath{\text{on}}\ensuremath{\ \Gamma_t}} 
\end{aligned}
\end{equation}
where $\normal$ is {the unit} normal vector to the interface $\Gamma_t$.
Here, $\rho$ and $\mu$ are positive functions that may have discontinuity across the interface, and $\left[q\right]=q^+-q^-$ denotes the jump relation, where the superscript $\pm$ refers to $\Omega^\pm$. In this paper, we {assume the following three conditions}: (\romannum{1}) $\left[\uvec\right]=a=0$, (\romannum{2}) interface $\Gamma_t$ moves with velocity $V$, and (\romannum{3}) $\rho$ and $\mu$ are piecewise constant.
{Our interest in moving interface problems \eqref{eq:PDE}--\eqref{eq:jmp} is motivated by a desire to develop a second-order accurate method for incompressible two-phase flows \cite{sussman1994level}. The time discretization to handle incompressibility and the pressure term for the two-phase flows are important and interesting topics; however, we cannot overlook spatial discretization with jump conditions. Therefore, we consider a simplified equation \eqref{eq:PDE} with jump conditions \eqref{eq:jmp}, which lacks the pressure term and incompressibility condition when compared with the two-phase flows. Nonetheless, it is very challenging to obtain a high-order accurate numerical method for \eqref{eq:PDE}--\eqref{eq:jmp} because it is sensitive to the determination of the moving interface, which is strongly coupled to the solution and the jump condition.} Thus, the aim of this study is to introduce a second-order accurate finite difference method for moving interface problems.

There are several numerical methods designed to treat jump conditions at fixed interfaces. One famous example is the immersed interface method (IIM) introduced by Leveque and Li \cite{leveque1994immersed}. The main idea of the IIM is to include the jump conditions in the discretization by utilizing a multivariate Taylor expansion. This was later extended in \cite{chen2019direct,li2017accurate,li1998fast}, thereby obtaining a second-order accurate solution and gradient at the interface. Recently, IIM with global second-order convergence in gradients was developed in \cite{tong2020obtain}. 
Another is the ghost fluid method (GFM) developed by Liu et al. in \cite{liu2000boundary}. The GFM introduces fictitious points to impose sharp jump conditions.
While the GFM is simpler to implement than IIM, it loses accuracy due to ignorance of tangential jump conditions. 
To address this issue, first- and second-order extensions of GFM were developed in \cite{egan2020xgfm} and \cite{guittet2015solving,cho2019second} for solving elliptic interface problems.
{Moreover,} there are other various methods for the interface problems on Cartesian grids. For examples, see the matched interface boundary method \cite{zhou2006high}, virtual node method \cite{bedrossian2010second,hellrung2012second}, ghost-point multi-grid method \cite{coco2018second}, and correction function method\cite{marques2011correction,marques2017high}.

These numerical approaches have been extended to handle non-smooth solutions across the moving interface. Li \cite{li1997immersed} proposed an IIM based numerical algorithm for solving the one-dimensional nonlinear moving interface problems. Second--order convergence was obtained for the solution and interface positions, and this work has been applied to the two--phase Navier--Stokes equations when density is constant in \cite{li2001immersed,lee2003immersed,tan2008immersed}. Some of these methods verify second-order convergence for non-smooth solutions. However, accuracy tests with the exact solution were only conducted with stationary interface; convergence with the moving interface was not reported. On the other hand, GFM motivated various sharp capturing methods for two-phase flows \cite{kang2000boundary,theillard2019sharp,sussman2007sharp,schroeder2014second}. These methods succeeded in capturing jump conditions for pressure and viscous terms; however, the jump conditions were not included in the discretization of convective terms. Some methods avoid this issue by extrapolating the velocities to another region. However, second-order convergence for piecewise-smooth solutions are only reported in \cite{schroeder2014second}.

In this paper, we introduce a second-order accurate semi-Lagrangian method for moving interface problems (\ref{eq:PDE})--\eqref{eq:jmp} on a Cartesian grid. The proposed method follows the framework of the second-order ghost fluid method \cite{cho2019second} to capture jump conditions and the level-set method \cite{osher1988fronts} to track the moving interface. The main idea of this study is to incorporate jump conditions into the discretization of convection terms and diffusion terms. Our main interest is the extension of this method to the development of a second-order accurate solver of two-phase flows. Thus, we propose a solution method and an interface tracking method for nonlinear systems of moving interface problems; that is, $V=\uvec$. The second-order convergence in both the solution and interface position is achieved.

The remainder of the paper is organized as follows: In section 2, we briefly review the level-set method and the second-order ghost fluid method. Section 3 contains the details of the numerical method for \eqref{eq:PDE}. Section 4 contains numerical experiments {that validate the second-order accuracy of the proposed method}.
\section{Preliminaries}

\subsection{Level-set Method}

In this study, the level-set method \cite{osher1988fronts} is used to represent and capture the moving interface. The interface $\Gamma_t$ is represented as a zero-level-set of a continuous function $\phi :  \Omega\times\left[0,T\right] \to \mathbb{R}$:
\begin{equation*}
\Gamma_t =  \left\{ \mathbf{x}\in \Omega  \mid  \phi( \mathbf{x},t ) =0 \right\}.
\end{equation*}
Furthermore, two subdomains $\Omega^+_t$ and $\Omega^-_t$ at time $t$ can be distinguished by the sign of the level-set function:
\begin{equation*}
\Omega^+_t =  \left\{ \mathbf{x}\in \Omega  \mid  \phi( \mathbf{x},t )  > 0 \right\}, \quad \Omega^-_t =  \left\{ \mathbf{x}\in \Omega \mid  \phi( \mathbf{x},t )  < 0 \right\}.
\end{equation*}
Evolution of the interface with the velocity ${V}$ is mathematically formulated as the following advection equation
\begin{equation} \label{eq:lvset}
\phi_{t}+{V}\cdot\nabla\phi=0.
\end{equation}

Another advantage of the level-set method is that the normal vector $\mathbf{n}$ can be simply represented in terms of the level-set function:
\begin{equation*}
\begin{aligned}
\normal&= (n_x, n_y) = \left( \frac{\phi_x}{\sqrt{\phi_x^2+\phi_y^2}},\frac{\phi_y}{\sqrt{\phi_x^2+\phi_y^2}}\right).
\end{aligned}
\end{equation*}
For more details, see \cite{gibou2018review,osher2004level}.
\subsection{\label{subsec:GFM2}A second-order ghost fluid method by Cho et al. \cite{cho2019second}}
In this subsection, we briefly review the second-order ghost fluid method proposed by Cho et al. \cite{cho2019second}, focusing on the core idea of approximating $u$ at the interface. Consider a jump condition
\[
\left[u\right]=0,	[\mu \frac{\partial \uvec}{\partial \normal}]=b   
\]
on an interface $\Gamma$.

\begin{figure}
	\centering{}\includegraphics[width=0.5\textwidth]{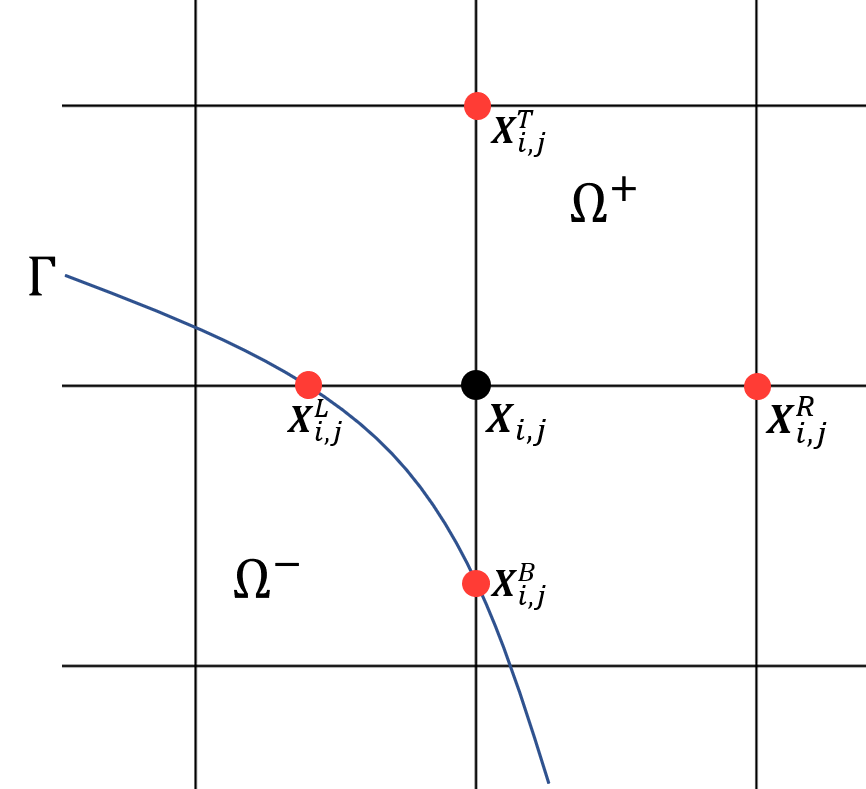}\caption{\label{fig:gfm2}Local four points.}
\end{figure}

Let $\pt_{ij}=\left(x_i,y_j\right)$ denote Cartesian grid points. Without loss of generality, assume $\pt_{ij} \in \Omega^+$. We define local five points $\ptij^{L},\ptij^{R},\ptij^{T}, \ptij^B$, and $\pt^\text{ext}_{ij}$ near $\pt_{ij}$. For example, $\ptij^R= \left(x_i+ \theta^R_{ij} \Delta x,y_j\right)$ for 
\[
\theta^{R}_{ij}=\begin{cases}
\frac{-D_{x}^{0}\phi-\text{sgn}\left(\phi_{i,j}\right)\sqrt{\left(D_{x}^{0}\phi\right)^{2}-4\phi_{i,j}D_{xx}^{0}\phi}}{2D_{xx}^{0}\phi} & \text{if }\phi_{i,j}\phi_{i+1,j}<0\\
1 & \text{if } \phi_{i,j}\phi_{i+1,j}>0
\end{cases}
\]
where $D_x^0\phi = \frac{\phi_{i+1,j}-\phi_{i-1,j}}{2}$ and $D_{xx}^0 \phi=\frac{\phi_{i+1,j}-2\phi_{i,j}+\phi_{i-1,j}}{2}$. Note that $\ptij^R= \pt_{i+1,j}$ if $\xy{i+1}{j}\in\Omega^+$. Conversely, if $\xy{i+1}{j}\in \Omega^-$, $\ptij^R$ is the interface location on a grid segment connecting $\pt_{ij}$ and $\xy{i+1}{j}$. As described in figure \ref{fig:gfm2}, 
\[\ptij^L= \left(x_i- \theta^{L}_{ij} \Delta x,y_j\right),\ptij^T= \left(x_i,y_j+ \theta^{T}_{ij} \Delta y\right),\ptij^B= \left(x_i,y_j- \theta^{B}_{ij} \Delta y\right)\]
are defined in a similar manner.

Let $\pt^\text{ext}_{ij}$ be any of the four points $\xy{i+l}{j+k}$ for $l,k=-1,1$ belonging to $\Omega^+$. Let $u_{ij},u^L_{ij},u^R_{ij},u^T_{ij},u^B_{ij}$ and $u^\text{ext}_{ij}$ be values of $u^+$ at six points $\pt_{ij},\ptij^L,\ptij^R,\ptij^T,\ptij^B$, and $\pt^\text{ext}_{ij}$, respectively. A quadratic polynomial $Q$ is constructed to interpolate the six values of $u^+$. Note that the coefficients of the quadratic polynomial $Q$ represent a linear combination of these six values of $u^+$.

For each $\ptij^L,\ptij^R,\ptij^T$ and $\ptij^B$, an equation is then set to approximate $u$ at the corresponding point. First, let us consider the equation corresponding to $\ptij^R$.
\begin{enumerate}
	\item If $\ptij^R=\pt_{i+1,j}$, simply equate $u^R_{ij}$ as $u$ at grid point $\pt_{i+1,j}$:
	\begin{equation}\label{eq:case1}
	u^R_{ij}=u_{i+1,j}.
	\end{equation}
	\item If $\ptij^R \neq \pt_{i+1,j}$, jump condition $\left[\mu u_x \right]$ is considered. Two different discretizations of $\left[\mu u_x \right]$ will be equated to construct the equation corresponding to $\ptij^R$. First, $\left[\mu u_x\right]$ can be expressed using the jump condition of normal derivatives:
	\begin{equation}
	[\mu u_x]= [\mu u_n]n_x-[\mu ]u_\tau^+  n_y,
	\end{equation}
	where $\tau$ is the tangent vector on the interface.
	By discretizing $u_\tau^+$ as $\nabla Q \cdot \tau$, one obtains a discretization of the jump condition $\left[\beta u_x \right]$. Another discretization of $\left[\beta u_x \right]$ is obtained a using one-sided second-order finite difference formula. Assume that $\pt_{i+1,j},\pt_{i+2,j}\in \Omega^-$, $u_x^+$, and $u_x^-$ at $\ptij^R$ are discretized using one-sided second-order finite difference formula:
	\begin{equation}\label{eq:second-order_finite_case2}
	\begin{aligned}
	u_x^+ &=\left.\left(\frac{\tR}{\tL(\tL+\tR)}u^{L}_{ij} -\frac{\tL+\tR}{\tL \tR}u_{i,j}+\frac{2\tR+\tL}{\tR(\tL+\tR)}u^{R}_{ij} \right)\middle/\Delta x \right. ,\\
	u_x^- &=\left.\left(-\frac{3-2\tR}{(1-\tR)(2-\tR)}u^{R}_{ij} +\frac{2-\tR}{1-\tR}u_{i+1,j}-\frac{1-\tR}{2-\tR}u_{i+2,j} \right)\middle/\Delta x \right. .
	\end{aligned}
	\end{equation}
	Two discretizations are combined to construct the equation corresponding to $\ptij^R$:
	\begin{equation}\label{eq:case2}
	[\beta u_n]n_x-[\beta ]\left( \nabla Q \cdot \tau \right)  n_y= \beta^+ u_x^+ - \beta^- u_x^-.
	\end{equation}
\end{enumerate}
In addition, an equation similar to either \eqref{eq:case1} or \eqref{eq:case2} is obtained for each $\ptij^L, \ptij^T$ and $\ptij^B$. Combining these four equations, we derive the following linear system 
\begin{equation}\label{eq:system_of_gfm2}
\mathbf{M} \begin{pmatrix}[1.5]
u^R_{ij}\\u^L_{ij}\\u^T_{ij}  \\ u^B_{ij} \end{pmatrix}
=\mathbf{N}\mathbf{u} + \mathbf{D},
\end{equation}
where $\mathbf{M}, \mathbf{N}, \mathbf{D}$ are real-valued vectors or matrices of appropriate size, and $\mathbf{u}$ is a vector consisting of values of $u$ at the grid points. By solving system \eqref{eq:system_of_gfm2}, we get expressions of $u^L_{ij},u^R_{ij},u^T_{ij},u^B_{ij}$ as linear combinations of $u$ at the grid points plus the constants. 

\paragraph{Remark} Note that it cannot be guaranteed that $\pt_{i+2,j}\in \Omega^-$ when discretizing \eqref{eq:second-order_finite_case2}. Assuming each domain $\Omega^+$ and $\Omega^-$ is connected, one simple approach to handle such a case is to set $\ptij^R=\xy{i+2}{j}$, $\theta_R=2$ and $u^R_{ij}= u_{i+2,j}$. For a detailed explanation of discretizing jump conditions when $\pt_{i+2,j}\in \Omega^+$, see \cite{cho2019second}.
\section{Numerical methods}
A collocated grid is used, so all values of $\phi$ and $u$ are located at Cartesian grid points $\pt_{ij}$, even when $u$ is vector-valued. To implicitly discretize the diffusive term in \eqref{eq:PDE}, the normal vector $\normal$ and the interface position at the next time level is needed. Thus, $\phi^{n+1}$ is updated before $u^{n+1}$. Once the level-set is advected, $u$ is solved using the semi-Lagrangian method combined with the backward difference formula. A concise outline of the overall algorithm is given as follows:	
\begin{enumerate}
	\item Evolve the level-set function and reinitialize it. (Before evolution of the level-set function, extrapolate $u^{n}$ at the interface to the grid points when $\mathbf{V}=u$)
	\item Trace back departure points and apply the interpolation procedure to discretize convective term via semi-Lagrangian method.
	\item Construct a linear system corresponding to backward difference discretization and solve it for $u$.
	
\end{enumerate} In the following sections, each step is explained in detail.

\subsection{Evolution of level-set}
The level-set function is updated from $\phi^{n}$ to $\phi^{n+1}$ with the velocity field $\vel$. To give the consistency with the overall methodology, we use the second-order semi-Lagrangian method
\begin{equation}\label{eq:lv_sl2}
\frac{3\phi_{ij}^{n+1}-4\phi_{d}^{n}+\phi_{d}^{n-1}}{2\Delta t}=0
\end{equation}
to discretize (\ref{eq:lvset}). The departure points are traced backward by a second-order Runge-Kutta method. The quadratic ENO interpolation procedure is then applied to recover the values of $\phi$ at the departure points. For more details, see \cite{min2007second}. 

The evolution (\ref{eq:lvset}) often leads numerical distortion of the level-set function near the interface. To avoid this problem, 
$\phi^{n+1}$ is reinitialized to the signed distance function, by solving the following pseudo-time dependent Eikonal equation:
\[
\phi_{\xi}+sgn\left(\phi^{0}\right)\left(\left\Vert \nabla\phi\right\Vert-1\right) =0.
\]
Here, $sgn$ denotes the signum function whose value is either -1, 0, or 1. For reinitialization, a Gauss-Seidel temporal discretization, in conjunction with ENO finite differences \cite{min2010reinitializing}, is used.

\subsubsection{Velocity extrapolation off the interface} \label{subsubsec:extrapolation}
In cases where $V=u$, the interface moves with velocity $u$. However, the discontinuity in gradient of $u$ causes first-order accuracy for the level-set function $\phi$, even with second-order semi-Lagrangian method and second-order reinitialization. Therefore, $\phi$ should be evolved with an alternative velocity field $W$ that has continuous gradients and agrees with $u$ on the interface, i.e. $W|_\Gamma=u|_\Gamma$. To construct such $W$, we adopt the extension algorithm technique with sub-cell resolution used in \cite{egan2020xgfm,dalmon2020fluids}.
The method extends $u$ to be a constant along the curve normal for $\Gamma$. 
This suggests the following pseudo-time dependent partial differential equation:

\begin{equation}\label{eq:psuedo_extra}
\begin{gathered}
\frac{\partial W}{\partial \xi} +\text{sign}(\phi)\normal \cdot \nabla W=0,\\
W=u \text{ on } \Gamma,
\end{gathered}
\end{equation}
whose characteristics are normal for $\Gamma$.

Assuming $\phi_{ij}>0$, the equation is semi-discretized as
\begin{equation*}
\frac{W^{l+1}_{ij} - W^l_{ij}}{\Delta \xi_{ij}} + (n_x^+ D_x^- W_{ij}+ n_x^- D_x^+ W_{ij}+n_y^+ D_y^- W_{ij}+ n_y^- D_y^+ W_{ij})=0,
\end{equation*}
where $n_x^+ = \max(n_x,0)$ and $n_x^- =\min(n_x,0)$.
Second-order ENO scheme with sub-cell resolution technique is used for spatial discretization.
\begin{equation*}
\begin{gathered}
D_x^- W_{ij} = \begin{cases}
\frac{W_{ij}-W_{i-1,j}}{\Delta x} + \frac{\text{minmod}\left( D_{xx}^0 W_{ij},D_{xx}^0 W_{i-1,j}\right)}{2\Delta x} &\text{if \ensuremath{\phi_{ij} \phi_{i-1,j} >0}}\\
\frac{W_{ij}-u^L_{ij}}{\theta^L_{ij} \Delta x} + \theta^L_{ij}\frac{\text{minmod}\left( D_{xx}^0 W_{ij},D_{xx}^0 W_{i-1,j}\right)}{2\Delta x}&\text{if \ensuremath{\phi_{ij} \phi_{i-1,j} \leq 0}}
\end{cases},\\
D_x^+ W_{ij} = \begin{cases}
\frac{W_{i+1,j}-W_{i,j}}{\Delta x} - \frac{\text{minmod}\left( D_{xx}^0 W_{ij},D_{xx}^0 W_{i-1,j}\right)}{2\Delta x}&\text{if \ensuremath{\phi_{ij} \phi_{i+1,j} >0}}\\
\frac{u^R_{ij}-W_{ij}}{\theta^R_{ij} \Delta x} -\theta^R_{ij} \frac{\text{minmod}\left( D_{xx}^0 W_{ij},D_{xx}^0 W_{i-1,j}\right)}{2\Delta x}&\text{if \ensuremath{\phi_{ij} \phi_{i+1,j} \leq 0}}
\end{cases}.
\end{gathered}
\end{equation*}
$D_y^- W_{ij}$ and $D_y^+ W_{ij}$ are similarly defined. We take a grid dependent time step $\Delta \xi_{ij}$ in order to avoid small time steps imposed at the grid points close to the interface. In particular, we set 
\begin{equation*}
\Delta \xi_{ij}= \min(\theta^R_{ij}, \theta^L_{ij}, \theta^T_{ij},\theta^B_{ij}) \times CFL \Delta x
\end{equation*}
And use the second-order Runge-Kutta formula for temporal discretization. In actual implementation, we take $CFL=0.4$ and the algorithm is performed up to 20 iterations.
After $W$ is obtained, the level-set function is advected with the extrapolated velocity $W$.

\subsection{\label{subsec:SL}Semi-Lagrangian ghost fluid method(SL-GFM)}
Quadratic and bilinear interpolation only attain first-order accuracy near the interface due to the discontinuity of gradients. Thus, we introduce a new interpolation scheme in conjunction with the ghost fluid method which provides second-order accuracy near the interface. 
Exploiting the jump conditions in computing the values of $u$ at
the departure points is the main idea for improving performance near
the interface. We name the semi-Lagrangian method with the proposed interpolation procedure the semi-Lagrangian ghost fluid method (SL-GFM).
\subsubsection{\label{subsubsec:RK2}Backward Integration}
We track the departure points by the second-order Runge-Kutta

\begin{align*}
\hat{\pt} & =\pt_{ij}-\frac{\Delta t}{2}\vel^{n}(\pt_{ij})\\
\pt_{d}^{n} & =\pt_{ij}-\Delta t\vel^{n+\frac{1}{2}}(\hat{\pt}).
\end{align*}
The velocity at the half time step $t^{n+\frac{1}{2}}$ is approximated
by the second-order extrapolation 
\[
\vel^{n+\frac{1}{2}}=\frac{3}{2}\vel^{n}-\frac{1}{2}\vel^{n-1}.
\]
Similarly, the departure point at the time step $n-1$ is estimated by
\begin{align*}
\hat{\pt} & =\pt_{ij}^{n+1}-\Delta t\vel^{n}(\pt_{ij})\\
\pt_{d}^{n-1} & =\pt_{ij}^{n+1}-2\Delta t\vel^{n}(\hat{\pt}).
\end{align*}
We impose a restriction on the time step size as 
\[
\Delta t<\frac{\Delta x}{2\parallel V\parallel_{\infty}}
\]
to ensure the departure points locate in an adjacent cell. That is, 
\begin{align*}
\mid x_{d}-x_{i}\mid &<\Delta x\\
\mid y_{d}-y_{j}\mid &<\Delta y.
\end{align*}

\subsubsection{Interpolation procedures} \label{subsubsec:interp_SL}

Let us focus on the interpolation procedure of the semi-Lagrangian method. Let $u_d^n$ denote the approximation of $u$ at the point $\pt_d^n$. Without loss of generality, assume $\pt_d^n\in\left[x_i, x_{i+1}\right]\times\left[y_j,y_{j+1}\right]$. Although $\pt_{ij}\in\Omega^+_{t^{n+1}}$ does not guarantee that $\pt_d^n\in \Omega^+_{t^{n}}$ and $\pt_d^{n-1}\in \Omega^+_{t^{n-1}}$ due to numerical error, we assume $\phi(\pt_d^n , t^n)>0$ if $\phi_{ij}^{n+1}>0$ and $\phi(\pt_d^n , t^n)<0$ if $\phi_{ij}^{n+1}<0$ during the interpolation procedure. For simplicity, in this section we omit subscripts and superscripts corresponding to time. The most natural choice for approximating $u_d$ is the bilinear interpolation 
\begin{equation}
u_{d} =u_{ij}\left(1-\theta_{x}\right)\left(1-\theta_{y}\right)+u_{i+1,j}\theta_{x}\left(1-\theta_{y}\right)+u_{i,j+1}\left(1-\theta_{x}\right)\theta_{y}+u_{i+1,j+1}\theta_{x}\theta_{y}, 
\label{eq:bilinear}
\end{equation}
where $\theta_{x}=\frac{x_{d}-x_{i}}{\Delta x},\ \theta_{y}=\frac{y_{d}-y_{j}}{\Delta y}\in\left(0,1\right)$.

We can construct a quadratic interpolation by correcting (\ref{eq:bilinear}) with second-order derivatives. Let us define a discrete operator 
\begin{equation}
D_{xx}^0 u_{ij}= u_{i+1,j}- 2u_{ij}+ u_{i-1,j}.
\end{equation}
For numerical stability of the interpolation, $u_{xx}^0$ is set to be one of $D_{xx}^0 u_{ij},D_{xx}^0 u_{i+1,j},D_{xx}^0 u_{i,j+1}$, and $D_{xx}^0 u_{i+1,j+1}$, which has the minimum absolute value; $u_{yy}^0$ is defined in duplication. These values lead to quadratic ENO interpolation
\begin{equation}
\begin{aligned} 
u_{d}  &= u_{ij}\left(1-\theta_{x}\right)\left(1-\theta_{y}\right)+u_{i+1,j}\theta_{x}\left(1-\theta_{y}\right)\\
&+ u_{i,j+1}\left(1-\theta_{x}\right)\theta_{y}+u_{i+1,j+1}\theta_{x}\theta_{y}-u_{xx}^0 \frac{\theta_{x}(1-\theta_{x})}{2}-u_{yy}^0 \frac{\theta_{y}(1-\theta_{y})}{2}.
\end{aligned} \label{eq:quadENO}
\end{equation}

The interpolation \eqref{eq:quadENO} is third-order accurate for smooth $u$. However, one cannot guarantee that all grid points involved in the interpolation belong to same subdomain with $\pt_d$. To address this issue, we introduce the following notions. We say the point $\pt_{d}$ is \it regular \rm if all grid points involved in the quadratic ENO interpolation belong to the same region with $\pt_{d}$. Otherwise, $\pt_d$ is called \it irregular \rm. In SL-GFM, the quadratic ENO interpolation is used for regular $\pt_d$ and the modified bilinear interpolation is used for irregular $\pt_d$. The modified bilinear interpolation is constructed by replacing the interpolating value to the ghost value as follows.

Suppose $\pt_d$ is irregular and there exists at least one grid point out of $\pt_{ij},\pt_{i+1,j},\pt_{i,j+1},\pt_{i+1,j+1}$ that belongs to the same region with $\pt_d$. Without loss of generality, we may assume $\pt_d\in\Omega^+$ and $\pt_{ij}\in \Omega^-$. Next, using $u_{ij}$ with the bilinear interpolation \eqref{eq:bilinear} would produce $O(\Delta x)$ error. Alternatively, we replace this with $u_{i,j}^+$, which is the extended value of $u^+$ at $\pt_{ij}$ obtained by using the second-order GFM. 

If $\xy{i+1}{j}\in\Omega^+$, the point $\ptij^R$ is located on the interface. One possible ghost value of $u^+_{ij}$ can be computed by an extrapolation from $u_{ij}^R$ and $u_{i+1,j}$:
\begin{equation*}
\frac{\theta^R_{ij}}{1-\theta^R_{ij}} (u_{ij}^R -u_{i+1,j}) +u_{ij}^R.
\end{equation*} Similarly, if $\pt_{i,j+1}\in \Omega^+$, $u^+_{ij}$ can be approximated as 	\begin{equation*}
\frac{\theta^T_{ij}}{1-\theta^T_{ij}} (u_{ij}^T -u_{i,j+1}) +u_{ij}^T.
\end{equation*} When the extrapolations from the both directions are available, as shwon in figure \ref{fig:grid_location}, the direction with the smaller distance is chosen. When both $\pt_{i+1,j},\pt_{i,j+1}\in \Omega^-$, as shown in figure \ref{fig:grid_location2}, $u^+_{ij}$ is approximated from $u_{i+1,j+1}$:
\begin{equation*}
u_{i+1,j+1} + \frac{1}{1-\theta^R_{i,j+1}}\left(u^R_{i,j+1}-u_{i+1,j+1}\right) +\frac{1}{1-\theta^T_{i+1,j}}\left(u^T_{i+1,j}-u_{i+1,j+1}\right).
\end{equation*}
In summary we derive the following formula:
\begin{figure} 
	\begin{center}
		
		\subfloat[][\label{fig:grid_location}]{\includegraphics[width=0.5 \textwidth]{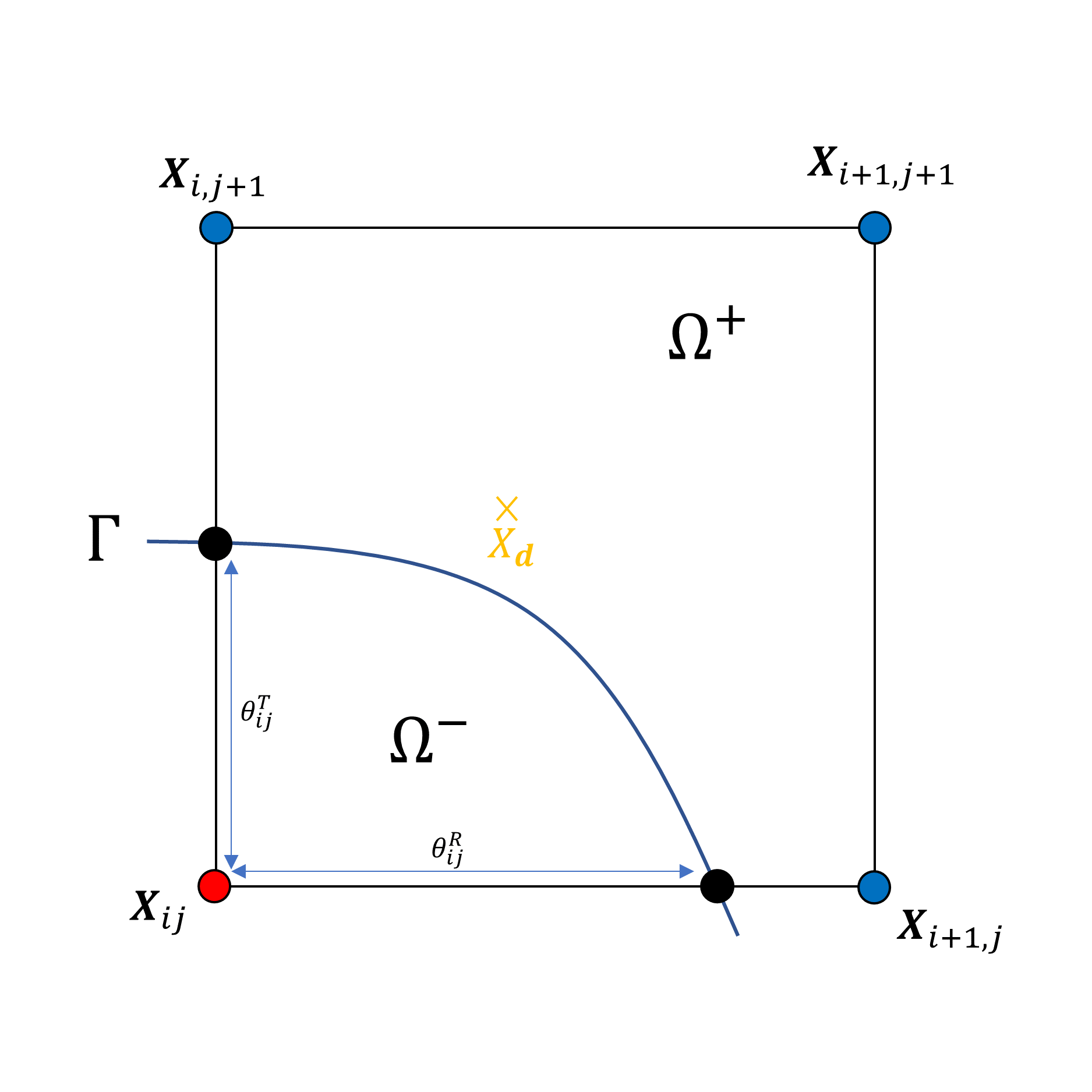}}
		\subfloat[][\label{fig:grid_location2}]{\includegraphics[width=0.5 \textwidth]{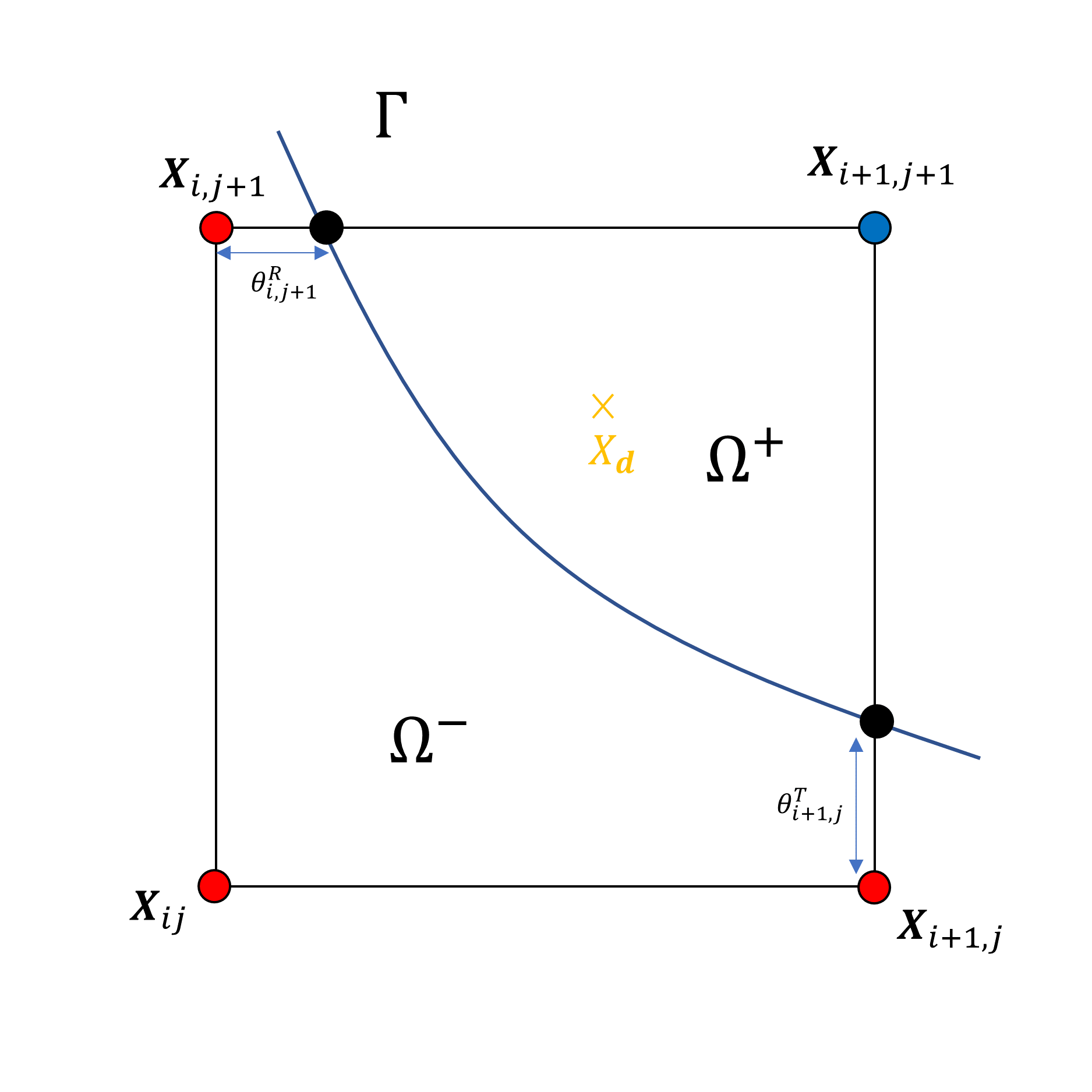}}
		\caption{ Grid points near the interface}
	\end{center}
\end{figure}
\[
u^+_{ij}=\begin{cases}
u_{ij} & \text{if \ensuremath{\ptij\in\Omega^+}}\\
\frac{1}{1-\theta^R_{ij}}u_{ij}^R - \frac{\theta^R_{ij}}{1-\theta^R_{ij}} u_{i+1,j} & \text{else if \ensuremath{\pt_{i+1,j} \in \Omega^+ \text{ and } \theta^R_{ij} \leq  \theta^T_{ij} }}\\
\frac{1}{1-\theta^T_{ij}}u_{ij}^T - \frac{\theta^T_{ij}}{1-\theta^T_{ij}} u_{i,j+1}  & \text{else if \ensuremath{\pt_{i,j+1}\in \Omega^+ \text{ and } \theta^T_{ij} \leq  \theta^R_{ij} }}\\
\frac{\theta^R_{i,j+1}\theta^T_{i+1,j}-1}{(1-\theta^R_{i,j+1})(1-\theta^T_{i+1,j})}	 u_{i+1,j+1} + \frac{1}{1-\theta^R_{i,j+1}}u^R_{i,j+1} +\frac{1}{1-\theta^T_{i+1,j}}u^T_{i+1,j} & \text{otherwise.}
\end{cases}.
\]
A similar process is used to define the ghost values $u^+_{i+1,j},u^+_{i,j+1}$
and $u^+_{i+1,j+1}$, which correspond to points $\xy{i+1}{j},\xy{i}{j+1},$
and $\xy{i+1}{j+1}$, respectively. Following this, we adopt the bilinear interpolation to approximate $u_d$ at $\Omega^+$:
\begin{equation}
u_{d}  =u^+_{ij}\left(1-\theta_{x}\right)\left(1-\theta_{y}\right)+u^+_{i+1,j}\theta_{x}\left(1-\theta_{y}\right)+u^+_{i,j+1}\left(1-\theta_{x}\right)\theta_{y}+u^+_{i+1,j+1}\theta_{x}\theta_{y}.
\label{eq:bilinear_mod}
\end{equation}

Due to numerical errors, it is not guaranteed that at least one of $\pt_{ij},\pt_{i+1,j},\pt_{i,j+1}, or \pt_{i+1,j+1}$ belongs to $\Omega^+$. In other words, the ghost values of $u^+$ may not be defined, so the interpolation from $\Omega^+$ is not possible on the cell. In such a case, $u_d$ is approximated via the bilinear interpolation on the four points in $\Omega^-$. To justify the approximation, we briefly prove the following statement:
\begin{equation*}
u^+(\pt_d)= u_d +O(\Delta x^2).	
\end{equation*}
Let 
\begin{equation*}
\phi_{d} =\phi_{ij}\left(1-\theta_{x}\right)\left(1-\theta_{y}\right)+\phi_{i+1,j}\theta_{x}\left(1-\theta_{y}\right)+\phi_{i,j+1}\left(1-\theta_{x}\right)\theta_{y}+\phi_{i+1,j+1}\theta_{x}\theta_{y}.
\end{equation*}
Since the bilinear interpolation is second-order accurate, we have 
\begin{equation}\label{eq:rem1}
| \phi_d - \phi(\pt_d)|=O(\Delta x^2).
\end{equation}
In addition, the condition $\pt_{ij},\pt_{i+1,j},\pt_{i,j+1},\pt_{i+1,j+1}\in \Omega^-$ leads to 
\begin{equation}\label{eq:rem2}
\phi_d <0.
\end{equation}
From (\ref{eq:rem1}), (\ref{eq:rem2}) together with $\pt_d \in \Omega^+$, we obtain 
\[
\phi(\pt_d)= O(\Delta x^2).
\]
Since $\phi$ is a signed distance function, there exist $\pt_\Gamma\in \Gamma$, such that $|\pt_d-\pt_\Gamma|=O(\Delta x^2)$. Thus, we may conclude that
\begin{equation*}
\begin{aligned}
|u_d - u^+(\pt_d)|  &\leq  |u_d -u^-(\pt_d)| +|u^-(\pt_d)-u^-(\pt_\Gamma)| +  |u^+(\pt_\Gamma)-u^+(\pt_d)|\\
&\leq O(\Delta x^2) + O(\Delta x^2) +O(\Delta x^2)=O(\Delta x^2).
\end{aligned}
\end{equation*}
Here, we used the facts that $|u(\pt_\Gamma)-u(\pt_d)| \leq |\nabla u |_{\infty} |\pt_\Gamma- \pt_d|+O(\Delta x^2)=O(\Delta x^2)$ and $u^-(\pt_\Gamma )=u^+(\pt_\Gamma)$.

\subsection{Linear system}
When $\pt_d^n$ and $\pt_d^{n-1}$ are both regular with respect to time $t^n$ and $t^{n-1}$, a second-order backward differentiation formula(BDF) is
adopted. Namely, 
\[
\rho \frac{3u_{ij}^{n+1}-4u_{d}^{n}+u_{d}^{n-1}}{2\Delta t}=\mu \bigtriangleup_h u^{n+1}_{ij}+f_{ij}^{n+1}.
\]
On the other hand, if one of $\pt_d^n$ and $\pt_d^{n-1}$ is irregular, first-order BDF is used for time discretization:
\[
\rho \frac{u_{ij}^{n+1}-u_{d}^{n}}{\Delta t}=\mu \bigtriangleup_h u^{n+1}_{ij}+f_{ij}^{n+1}.
\]
Following the second-order GFM, the standard finite difference method of the five-point Laplacian formula is used to discretize the elliptic operator at the grid points away from the interface $\Gamma_{t^{n+1}}$:
\begin{equation*}
\begin{aligned}
\bigtriangleup_h u_{ij}^{n+1} = \frac{ u^{n+1}_{i+1,j}-2u^{n+1}_{ij}+u^{n+1}_{i-1,j} }{\Delta x^2} + \frac{ u^{n+1}_{i,j+1}-2u^{n+1}_{ij}+u^{n+1}_{i,j-1} }{\Delta y^2}.
\end{aligned}
\end{equation*}	
If grid segments for each Cartesian direction intersect with the interface, the Laplacian formula is discretized using the Shortley-Weller method \cite{shortley1938numerical}
\begin{equation*}
\begin{aligned}
\bigtriangleup_h u_{ij}^{n+1} &= \left. \left(\frac{ u^{R}_{ij}-u^{n+1}_{ij}}{\theta^R_{ij} \Delta x} -\frac{u_{ij}^{n+1}-u^{L}_{ij} }{\theta^L_{ij} \Delta x} \right)\middle/\left(\frac{\theta^R_{ij}+\theta^L_{ij}}{2}\Delta x\right) \right.\\
&+ \left.\left(\frac{ u^{T}_{ij}-u^{n+1}_{ij}}{\theta^T_{ij} \Delta x} -\frac{u_{ij}^{n+1}-u^{B}_{ij} }{\theta^B_{ij} \Delta y } \right)\middle/\left(\frac{\theta^T_{ij}+\theta^B_{ij}}{2}\Delta y\right),\right.
\end{aligned}
\end{equation*}
where $u_{ij}^R,u_{ij}^L,u_{ij}^T$ and $u_{ij}^B$ are computed according to section \ref{subsec:GFM2}.
$\mu$ and $\rho$ are defined naturally, as 
\begin{equation*}
\rho = \begin{cases}
\rho^+ & \text{if \ensuremath{\ptij\in\Omega^+_{t^{n+1}}}}\\
\rho^- &\text{if \ensuremath{\ptij\in\Omega^-_{t^{n+1}}}}
\end{cases}, \quad
\mu =\begin{cases}
\mu^+ & \text{if \ensuremath{\ptij\in\Omega^+_{t^{n+1}}}}\\
\mu^- &\text{if \ensuremath{\ptij\in\Omega^-_{t^{n+1}}}}
\end{cases}.
\end{equation*}
Combining all, we attain a linear system whose coefficient matrix is non-symmetric. Since $u^R_{ij},u^L_{ij},u^T_{ij}$ and $u^B_{ij}$ are used in the interpolation at the next time step, we save these in actual implementation.
\paragraph{Remark}
The numerical method proposed throughout previous sections allows us to solve parabolic moving interface problems
\[
\rho u_t-\nabla\cdot\left(\mu \nabla u\right)=f.
\]
Except for the semi-Lagrangian application owing to the absence of advection terms, the second-order BDF
\[
\rho\frac{3u_{ij}^{n+1}-4u_{ij}^{n}+u_{ij}^{n-1}}{2\Delta t}=\nabla\cdot\left(\beta\nabla u_{i}^{n+1}\right)+f_{i}^{n+1}
\] 
can be used to design a second-order accurate scheme. If $\phi_{ij}^{n+1}>0$ but $\phi_{ij}^n<0$, implicit time discretization with ghost value $u^+_{ij}$, which was introduced in \ref{subsubsec:interp_SL}, is used:
\[
\rho \frac{u_{ij}^{n+1}-{u}_{ij}^{+}}{\Delta t}=\nabla\cdot\left(\beta\nabla u_{i}^{n+1}\right)+f_{i}^{n+1}.
\]

\section{Numerical results}

In this section, several numerical experiments are carried out to verify
second-order accuracy in the $L^{\infty}$ norm of the proposed method. 
Throughout this section, the linear system in section 3.3 is solved by the generalized minimal residual method (GMRES) with an incomplete LU preconditioner \cite{saad2003iterative}.
All numerical experiments are carried out in C++ on a personal computer. Unless the source term $f$ or the jump condition $\left[\mu \frac{\partial u}{\partial \normal}\right]$ is specifically mentioned, it is computed according to the exact solution $u$ and $\phi$.

\subsection{Scalar equation : Translation}
\begin{figure}
	\centering{}
	\subfloat[width=0.5\textwidth][Graph of $u$ by SL-GFM]{\includegraphics[width=0.5\textwidth]{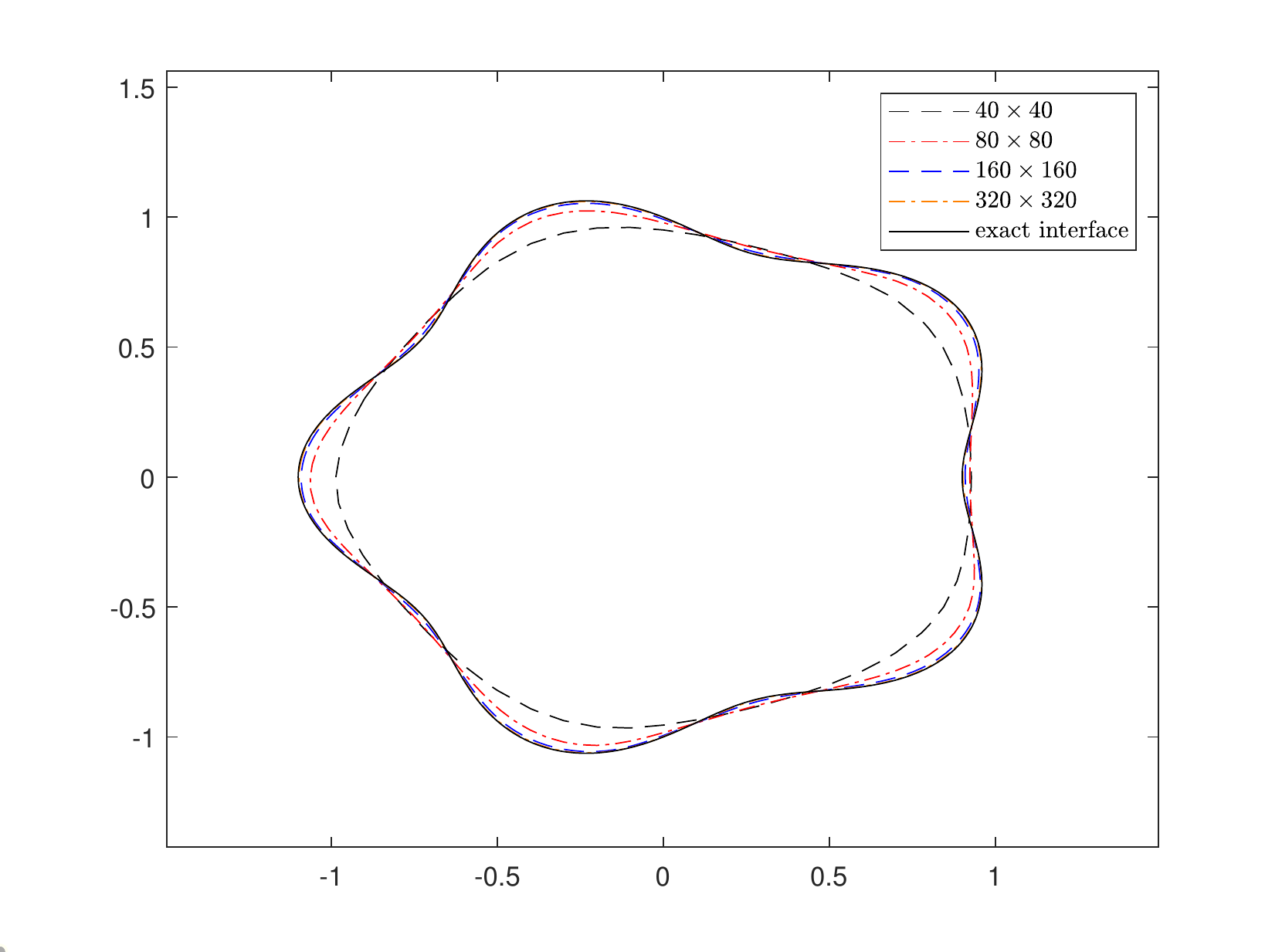}}
	\subfloat[width=0.5\textwidth][Cross section of $u$ at $y=1$ by SL-GFM(diamond) and SL-BDF2(triangle) with the exact value(solid line)]{\includegraphics[width=0.5\textwidth]{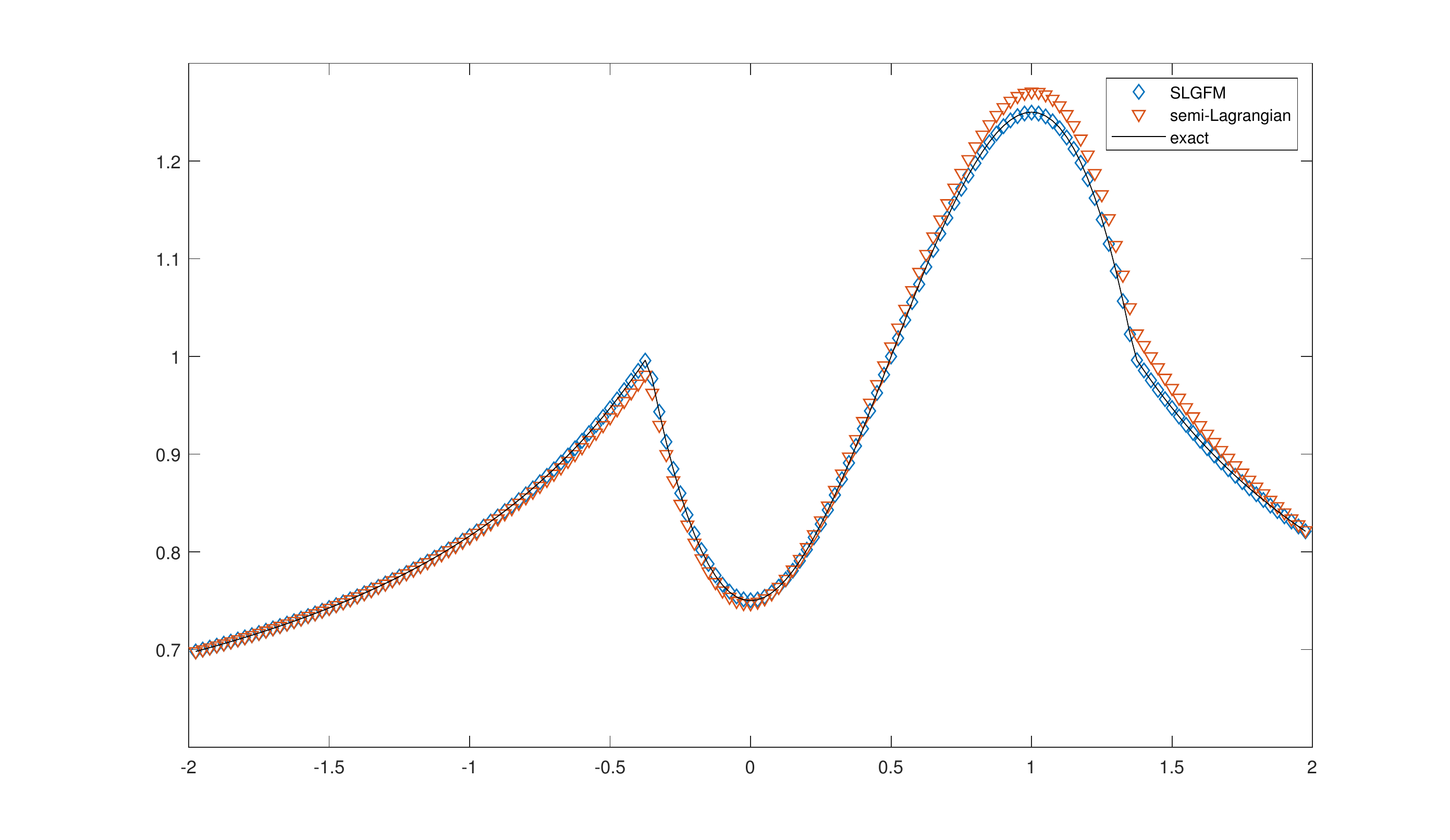}}
	\caption{\label{fig:result_circle}Solution of test 1 on a $160\times160$ grid.}
\end{figure}

We begin with an accuracy test in a simple setting. 
Consider a translating circular interface $\Gamma_t$ with velocity ${V}=\left(1,1\right)$ on a computational domain $\Omega=\left[-2,2\right]^2$. Specifically, the interface $\Gamma_t$ is defined as the zero level-set of a function 
\[
\phi(x,y,t)=\sqrt{\left(x-t+0.5\right)^{2}+\left(y-t+0.5\right)^{2}}-1.
\]
With the parameters
\begin{equation*}
\rho^-=1, \quad\rho^+=1, \quad\mu^-=1, \quad\mu^+=2
\end{equation*}
and the source term
\begin{equation*}
f(x,y,t)=\begin{cases}
4 \mu^-(2x - 2 t + 1);  & \text{in }\Omega^{-}\\
\frac{\mu^+( 2 y-2t+1)}{2((x-t+0.5)^2 +(y-t+0.5)^2)^\frac{3}{2}} & \text{in }\Omega^{+}
\end{cases},
\end{equation*}
the exact solution is given by
\[
u(x,y,t)=\begin{cases}
-(x-t+0.5) \left( (x-t+0.5)^2+(y-t+0.5)^2 - 1 \right))  & \text{in }\Omega^{-}\\
\frac{ y-t+0.5 }{\sqrt{(x-t+0.5)^2 + (y-t+0.5)^2}} - (y-t+0.5)  & \text{in }\Omega^{+}
\end{cases}.
\]
We compare the results obtained by our method and SL-BDF2, which uses the conventional quadratic ENO interpolation
with the semi-Lagrangian method.
Table \ref{tab:2_1} shows the convergence results with the time step $\Delta t=0.4\Delta x$ and at the final time $T=1$. Figure \ref{fig:result_circle} shows the corresponding solutions on the $160\times160$ grid.
Though the SL-BDF2 also discretizes diffusive terms using a ghost-fluid method \cite{cho2019second}, overall accuracy is first-order.
Besides, we can see that the SL-GFM yields second-order convergence.
The difference between these two results is the applied interpolation scheme; hence, it indicates the necessity of using the proposed interpolation method.
It can also be confirmed from figure \ref{fig:result_circle} that the SL-GFM approximates the solution $u$ much more accurately than the SL-BDF2.

\begin{table}[H]
	\begin{centering}
		
		\begin{tabular}{|cccccc|}
			\hline 
			Grid & \multicolumn{2}{c}{SL-BDF2} &  & \multicolumn{2}{c|}{SL-GFM}\tabularnewline
			\cline{2-6} 
			& $\parallel u(x,y,t)-u_{ij}\parallel_{\infty}$ & order &  & $\parallel u(x,y,t)-u_{ij}\parallel_{\infty}$ & order\tabularnewline
			\hline 
			$40^2$                 & 1.34$\times 10^{-1}$ & - &  & 1.05$\times 10^{-2}$ & - \\ \hline
			$80^2$                 & 6.70$\times 10^{-2}$ & 1.00 &  & 2.62$\times 10^{-3}$ & 2.00 \\\hline
			$160^2$                & 3.29$\times 10^{-2}$ & 1.03 &  & 7.02$\times 10^{-4}$ & 1.90 \\\hline
			$320^2$                & 1.65$\times 10^{-2}$ & 1.00 &  & 1.82$\times 10^{-4}$ & 1.95 \\\hline
			$640^2$                & 8.27$\times 10^{-3}$ & 0.99 &  & 4.60$\times 10^{-5}$ & 1.98 \\
			\hline 
		\end{tabular}
		\par\end{centering}
	\caption{\label{tab:2_1}Convergence rates for test 1.}
\end{table}

\subsection{Scalar equation : Rotation}
As a second example, a flower shaped interface rotating around the origin with a unit angular velocity $V=\left(-y,x\right)$ in a domain $\Omega=[-2,2]^{2}$ is chosen.
The interface can be described by the level function in polar coordinates
\[
\phi(x,y,t)=r-1-0.1\cos\left(5\left(\theta-t\right)\right).
\]
Parameters are chosen to $\rho^{-}=100,\rho^{+}=1, \mu^{-}=10$ and $\mu^{+}=1$  with the
exact solution 
\[
u(x,y,t)=\begin{cases}
r^{2}-1 & \text{in }\Omega^{-}\\
0.1\cos(5(\theta-t))\left(2+0.1\cos(5(\theta-t))\right) & \text{in }\Omega^{+}.
\end{cases}
\]

\begin{figure}
	\centering{}
	\subfloat[width=0.5\textwidth][Contour of zero level-sets]{\includegraphics[width=0.5\textwidth]{./ex2_interf1}}
	\subfloat[width=0.5\textwidth][zommed]{\includegraphics[width=0.5\textwidth]{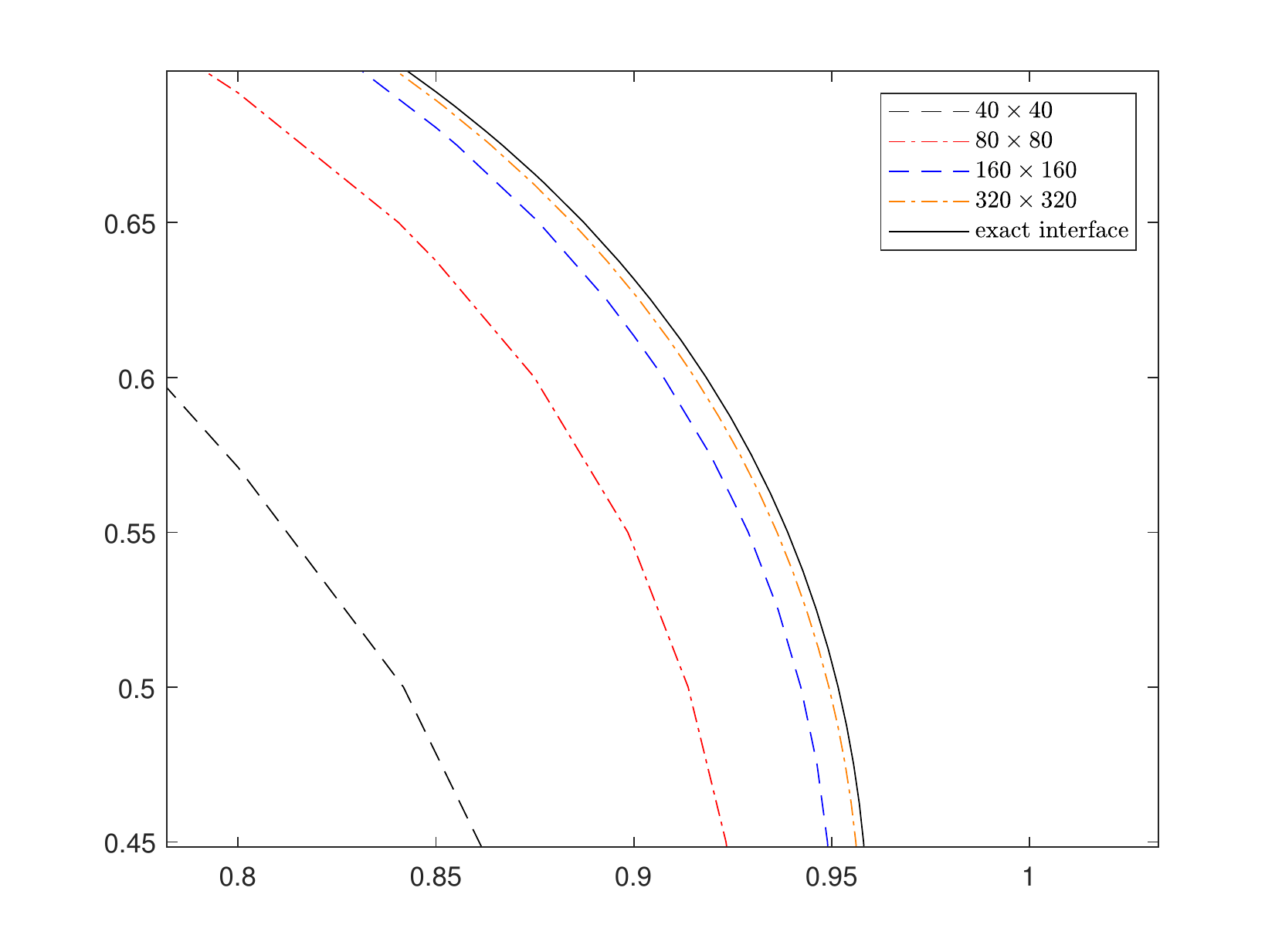}}
	\caption{\label{fig:result_flower}Interface position of example 2 at $t=\pi$ on various grids.}
\end{figure}

Since $\parallel V\parallel_{\infty}=4$ on $\Omega$, we set the time step restriction as $\Delta t=0.2 \Delta x$.
Results at the final time $T=\pi$, when the interface $\Gamma_t$ completes half a rotation,
are presented in table \ref{tab:2_2}. 
These results show that SL-GFM offers better accuracy when compared with SL-BDF2. The accuracy at low resolutions tends to be lower than second-order; however, this is due to the slow convergence of the interface position and the normal vector, which is evaluated from the level-set function. 
Referring to figure \ref{fig:result_flower}, convergence of the interface position is followed by the second-order convergence of the solution on finer resolution.                

\begin{table}[H]
	\begin{centering}	
		\begin{tabular}{|cccccc|}
			\hline 
			Grid & \multicolumn{2}{c}{SL-BDF2} &  & \multicolumn{2}{c|}{SL-GFM}\tabularnewline
			\cline{2-6} 
			& $\parallel u(x,y,t)-u_{ij}\parallel_{\infty}$ & order &  & $\parallel u(x,y,t)-u_{ij}\parallel_{\infty}$ & order\tabularnewline
			\hline 
			$40^{2}$ & $2.17\times10^{-1}$ & - &  & $1.62\times10^{-1}$ & -\tabularnewline
			\hline 
			$80^{2}$ & $1.34\times10^{-1}$ & 0.69 &  & $7.48\times10^{-2}$ & 1.11\tabularnewline
			\hline 
			$160^{2}$ & $6.61\times10^{-2}$ & 1.02 &  & $2.15\times10^{-2}$ & 1.80\tabularnewline
			\hline 
			$320^{2}$ & $2.95\times10^{-2}$ & 1.16 &  & $6.01\times10^{-3}$ & 1.84\tabularnewline
			\hline 
			$640^{2}$ & $1.42\times10^{-2}$ & 1.06 &  & $1.57\times10^{-3}$ &1.94 \tabularnewline
			\hline 
		\end{tabular}
		\par\end{centering}
	\caption{\label{tab:2_2}Convergence rates for test 2.}
\end{table}

\subsection{Scalar equation : Deformation}
\begin{figure}
	\centering{}\includegraphics[width=0.5\textwidth]{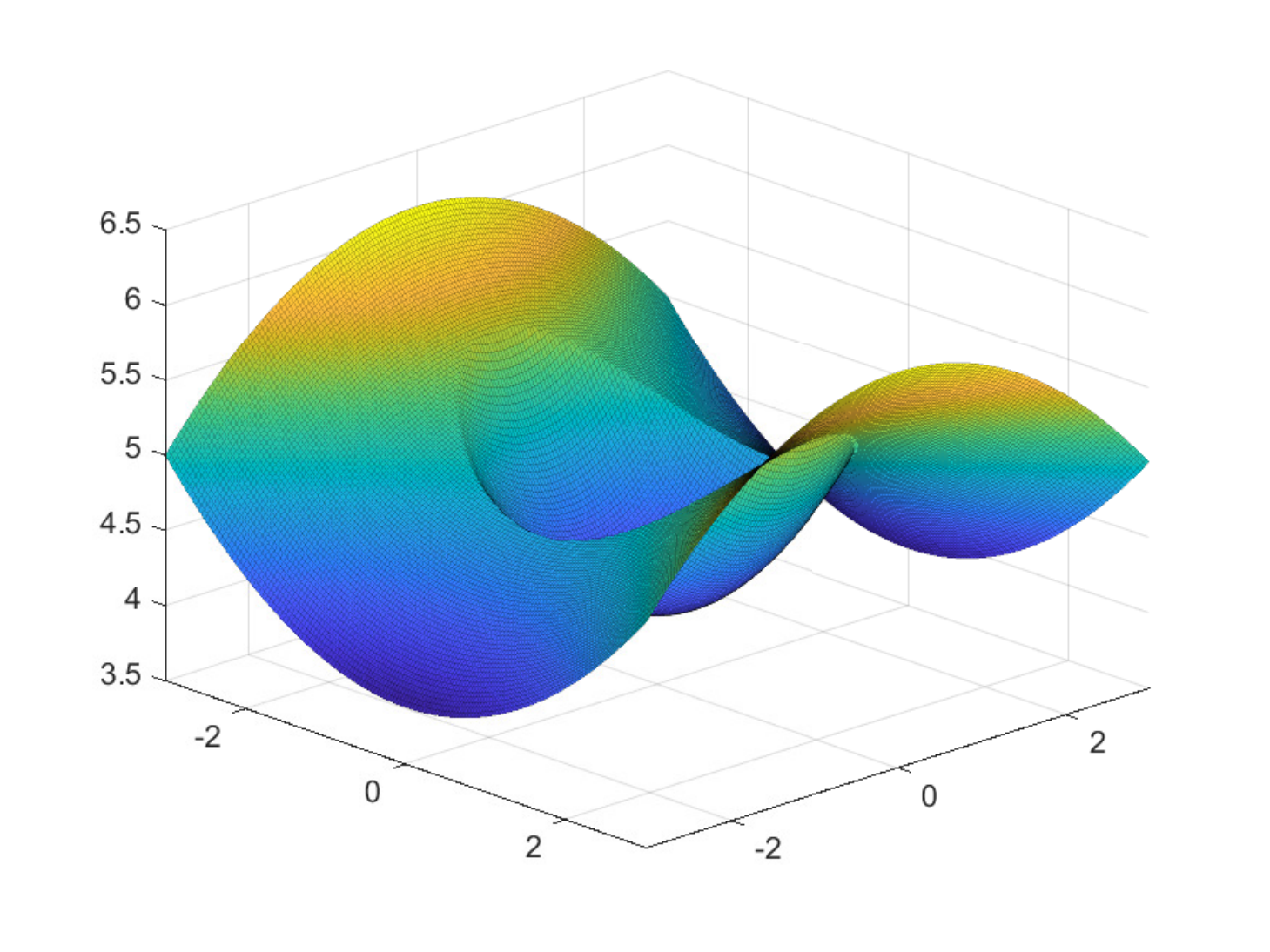}\caption{\label{fig:result_ellipse}Solution of test 3 on $160\times160$ grid.}
\end{figure}

Consider a deforming interface 
\begin{equation}
\Gamma_t =\{ \left(x,y \right) \in \mathbb{R}^2 |x^{2}\left(1-\frac{3}{4}t\right)+y^{2}\left(1-\frac{1}{2}t\right)=1 \} .
\end{equation} 
It is easy to see that $\Gamma_t$ is a unit circle at $t=0$ and moves with the velocity
\[
V(x,y,t)=\left[\frac{3x}{8-6t},\ \frac{y}{4-2t}\right],
\]
to become an ellipse. Two methods are tested with solution
\[
u(x,y,t)=\begin{cases}
\left(1-\frac{5}{8}t\right)\left(x^2+y^2\right)+4t & \text{in }\Omega^{-}\\
\frac{1}{8}t\left(x^2-y^2\right)+4t & \text{in }\Omega^{+}
\end{cases}.
\]
and quantities $\rho^{-}=1,\rho^{+}=100,\mu^{-}=1, \mu^{+}=10$.
In this simulation, we used the time step restriction $\Delta t=0.1\Delta x$.
The accuracy at the final time $T=1$ is presented in table \ref{tab:2_3}.
Figure \ref{fig:result_ellipse} depicts the numerical solution obtained by the SL-GFM carried out on a $160\times160$ grid. 
It shows that the SL-GFM leads to second-order convergence, but the SL-BDF2 only results in the first-order accuracy.

\begin{table}[H]
	\begin{centering}
		\begin{tabular}{|cccccc|}
			\hline 
			Grid & \multicolumn{2}{c}{SL-BDF2} &  & \multicolumn{2}{c|}{SL-GFM}\tabularnewline
			\cline{2-6} 
			& $\parallel u(x,y,t)-u_{ij}\parallel_{\infty}$ & order &  & $\parallel u(x,y,t)-u_{ij}\parallel_{\infty}$ & order\tabularnewline
			\hline 
			$40^{2}$ & $1.01\times10^{-1}$ & - &  & $2.72\times10^{-2}$ & -\tabularnewline
			\hline 
			$80^{2}$ & $7.77\times10^{-2} $ & 0.38 &  & $7.93\times10^{-3}$ & 1.78\tabularnewline
			\hline 
			$160^{2}$ & $4.59\times10^{-2} $ & 0.76 &  & $2.18\times10^{-3}$ & 1.86\tabularnewline
			\hline 
			$320^{2}$ & $2.60\times10^{-2}$ & 0.82 &  & $5.84\times10^{-4}$ & 1.90\tabularnewline
			\hline 
			$640^{2}$ & $1.43\times10^{-2}$ & 0.86 &  & $1.48\times10^{-4}$ & 1.97\tabularnewline
			\hline 
		\end{tabular}
		\par\end{centering}
	\caption{\label{tab:2_3}Convergence rates for test 3.}
\end{table}

\subsection{Non-linear system : Translation}
\begin{figure}
	\subfloat[$u^1$]{\includegraphics[width=0.5\textwidth]{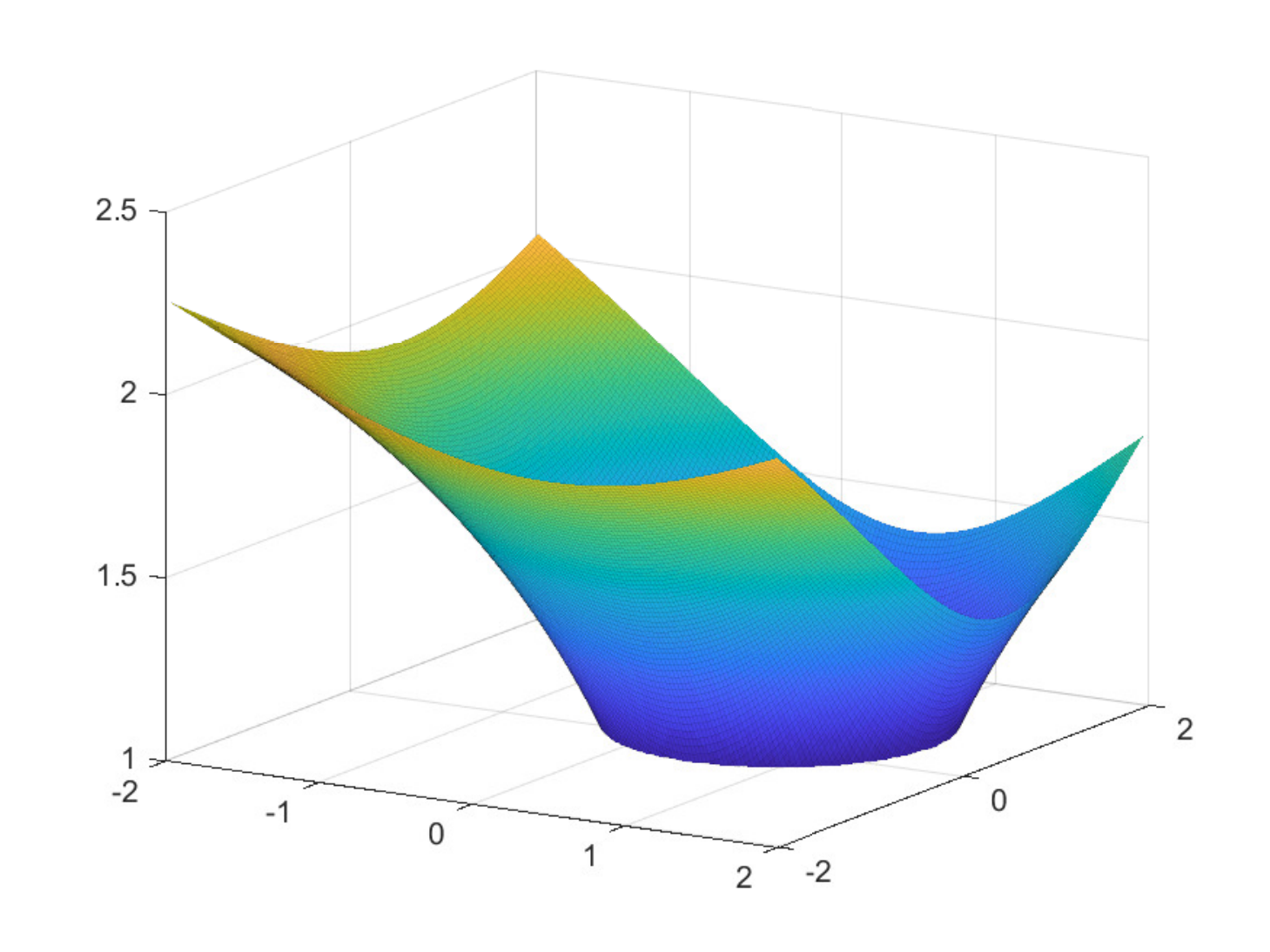}}
	\subfloat[$u^2$]{\includegraphics[width=0.5\textwidth]{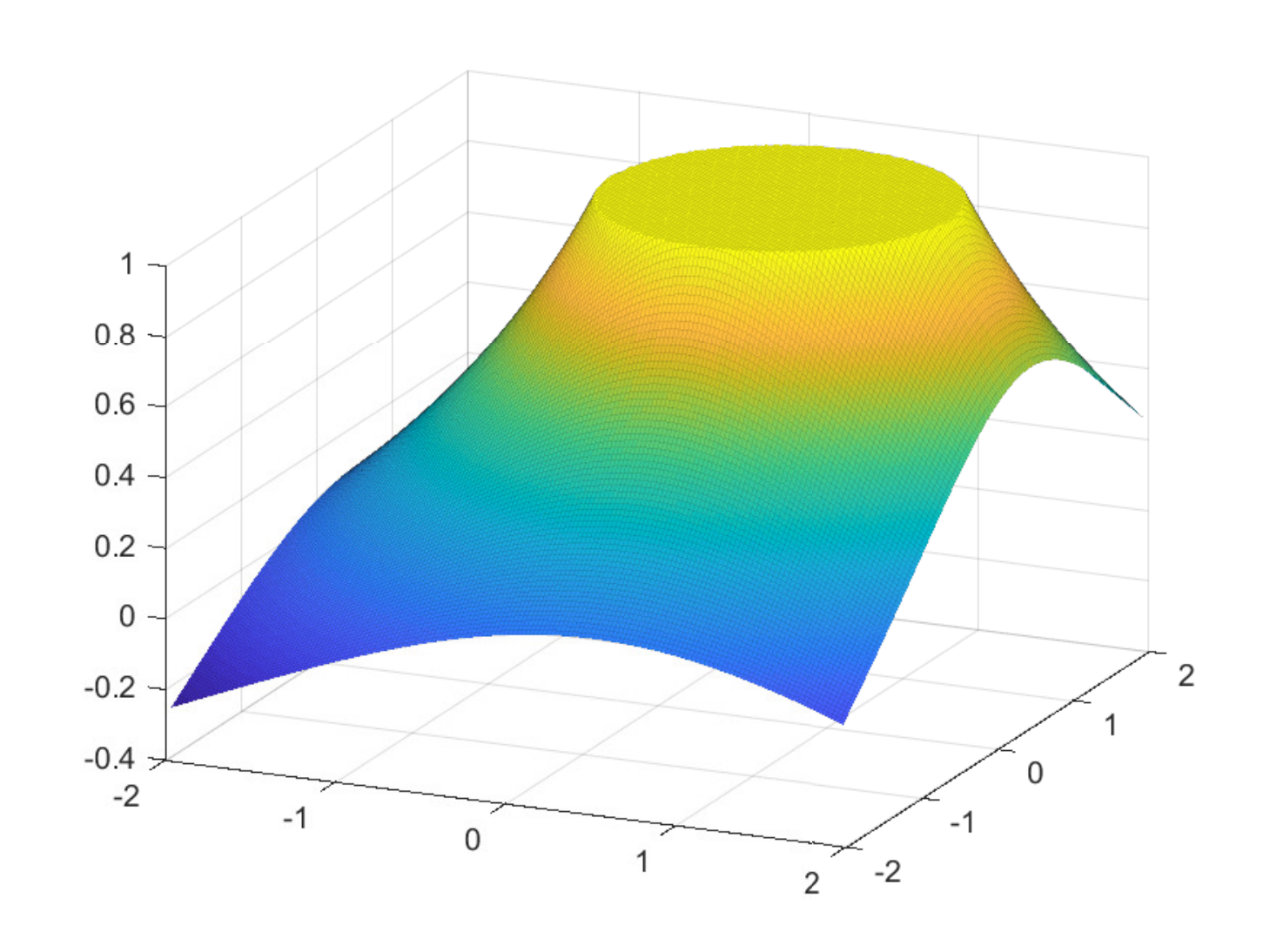}}
	\caption{\label{fig:ex4}Solution of test 4 on $160\times160$ grid.}
\end{figure}

We now consider a non-linear system, which occurs when $V=u$. On a computational domain $\Omega= \left[-2,2\right]^2$, the interface is given as the zero level-set of a function
\[
\phi(x,y,t)=\sqrt{\left(x-t+0.5\right)^{2}+\left(y-t+0.5\right)^{2}}-1,
\]
and the solution $u=\left(u^1, u^2\right)$ is given as 
\[
u^1(x,y,t)=\begin{cases}
1 & \text{in }\Omega^{-}\\
1+\frac{1}{2}\log \left( (x-t+0.5)^2 + (y-t+0.5)^2 )\right) & \text{in }\Omega^{+}
\end{cases}
\]
and
\[
u^2(x,y,t)=\begin{cases}
1 & \text{in }\Omega^{-}\\
1-\frac{1}{2}\log \left( (x-t+0.5)^2 + (y-t+0.5)^2 )\right) & \text{in }\Omega^{+}
\end{cases}.
\]
It is easy to see that $u|_\Gamma = \left(1,1\right)$, which agrees with the velocity of the interface. Jump conditions are
\[
\left[\mu \frac{\partial u}{\partial \normal} \right]= \left(\mu^+, \mu^- \right).
\]
Three different numerical experiments, using different methods for level-set advection, are conducted with $\rho^+=1,\rho^-=1000$ and $\mu^+=0.1, \mu^- =10$ up to $t=1$. See figure \ref{fig:ex4} for the profile of the solution.

First, let $u^a_{ij},\phi^a_{ij}$ denote numerical solutions of $u$ and $\phi$ when the level-set function is advected with the velocity $u^{n}$ at each step. Next, let $u^b_{ij},\phi^b_{ij}$ be numerical solutions when the extrapolation technique discussed in section \ref{subsubsec:extrapolation} is used to determine the velocity at which the level-set moves. Finally, $u^c_{ij}$ denotes a numerical solution of $u$ when $\phi$ is given exactly. Numerical errors for these solutions are presented in table \ref{tab:ex4}. In addition, $L^\infty$ errors of $\phi$ are computed only at the grid points, such that $|\phi|<3 \max(\Delta x, \Delta y)$. Despite the same method being applied to $u$ in these three tests, we can observe a huge difference between them due to the accuracy of tracking the interface. When $\phi$ is advected with $u$ but without the extrapolation technique, both $u^a_{ij}$ and $\phi^a_{ij}$ show first-order accuracy. However, when $\phi$ is computed with the extrapolated velocity, second-order convergence is obtained for both $u^b_{ij}$ and $\phi^b_{ij}$. Furthermore, when the interface is given exactly, we see that $u^c_{ij}$ shows the lowest error among all three experiments. Since the jump conditions $\left[\mu \frac{\partial u}{\partial \normal} \right]$ are dependent on normal vectors, a second-order accurate interface position and normal vector are needed to obtain a second-order accurate solution $u$.
\begin{table}[H]
	
	\centering
	\caption{Convergence rates for example 4\label{tab:ex4}}
	\subfloat[(a)][solution]{
		\begin{tabular}{|ccccccc|}
			\hline 
			Grid & $\parallel u(x,y,t)-u^a_{ij}\parallel_{\infty}$ & order &  $\parallel u(x,y,t)-u^b_{ij}\parallel_{\infty}$ & order& $\parallel u(x,y,t)-u^c_{ij}\parallel_{\infty}$ & order\tabularnewline
			\hline 
			$40^2$                 & $2.30\times10^{-2}$ & - & $1.24\times10^{-2}$ & - & $1.76\times10^{-3}$ & - \\ \hline
			$80^2$                 & $9.52\times10^{-3}$ & 1.27 & $3.19\times10^{-3}$ & 1.95 & $4.87\times10^{-4}$ & 1.85 \\\hline
			$160^2$                & $4.57\times10^{-3}$ & 1.06 & $9.13\times10^{-4}$ & 1.81 & $1.28\times10^{-4}$ & 1.92 \\\hline
			$320^2$                & $2.19\times10^{-3}$ & 1.06 & $2.53\times10^{-4}$ & 1.85 & $3.29\times10^{-5}$ & 1.97 \\\hline
			$640^2$                & $1.10\times10^{-3}$ & 1.00 & $6.33\times10^{-5}$ & 2.00 & $8.31\times10^{-6}$ & 1.98 \\
			\hline	
		\end{tabular}
	}
	
	\centering
	
	\subfloat[(b)][interface]{
		\begin{tabular}{|cccccc|}
			\hline  
			
			Grid & $\parallel \phi(x,y,t)-\phi^a_{ij}\parallel_{\infty}$ & order &  & $\parallel \phi(x,y,t)-\phi^b_{ij}\parallel_{\infty}$ & order\tabularnewline
			\hline 
			$40^2$                 & $2.37\times10^{-2}$ & - &  & $1.36\times10^{-2}$ & - \\ \hline
			$80^2$                 & $9.90\times10^{-3}$ & 1.26 &  & $3.33\times10^{-3}$ & 2.03 \\\hline
			$160^2$                & $4.88\times10^{-3}$ & 1.02 &  & $9.38\times10^{-4}$ & 1.83 \\\hline
			$320^2$                & $2.27\times10^{-3}$ & 1.11 &  & $2.56\times10^{-4}$ & 1.87 \\\hline
			$640^2$                & $1.14\times10^{-3}$ & 1.00 &  & $6.39\times10^{-5}$ & 2.00 \\\hline
			
		\end{tabular}
		
	}
	
\end{table}

\subsection{Non-linear system : Rotation}
Consider a circular interface $\Gamma_t$ that rotates around the origin, which is given as a zero level-set of a function
\[
\phi\left(x,y,t\right) = \sqrt{ (x-1.5\cos(t))^2 + (y-1.5\sin(t))^2} -1
\]
on computational domain $\left[-3,3\right]^2$. Hence, we conducted the accuracy test for the following solution
\[
u^1(x,y,t)=\begin{cases}
-y & \text{in }\Omega^{-}\\
-y+\frac{1}{2}\log \left( (x-1.5\cos(t))^2 + (y-1.5\sin(t))^2 )\right)& \text{in }\Omega^{+}
\end{cases}.
\]
and
\[
u^2(x,y,t)=\begin{cases}
x & \text{in }\Omega^{-}\\
x+\frac{1}{2}\log \left( (x-\cos(t))^2 + (y-\sin(t))^2 )\right) & \text{in }\Omega^{+}
\end{cases}.
\]
$f$ and $\left[\mu \frac{\partial u}{\partial \normal}\right]$ are given according to the solution.
Numerical simulations are performed up to $t= \pi$ with quantities $\rho^+ =1000,\rho^-=1, \mu^+=10, \mu^-=0.1$. Second-order convergences in $L^\infty$ norm of the solution and the interface position for both cases are presented in table \ref{tab:ex5}. 
\begin{table}[H]
	
	\centering
	\caption{Convergence rates for example 5\label{tab:ex5}}
	
	\begin{tabular}{|cccccc|}
		\hline 
		Grid & $\parallel u(x,y,t)-u_{ij}\parallel_{\infty}$ & order &  & $\parallel \phi(x,y,t)-\phi_{ij}\parallel_{\infty}$ & order\tabularnewline
		\hline 
		$40^2$   & 4.56$\times10^{-1}$ & - &  & 4.05$\times10^{-1}$ & - \\ \hline
		$80^2$   & 9.64$\times10^{-2}$ & 2.24 &  & 1.10$\times10^{-1}$ & 1.88 \\\hline
		$160^2$ & 2.80$\times10^{-2}$ & 1.78 &  & 3.03$\times10^{-2}$ & 1.86 \\\hline
		$320^2$ & 6.65$\times10^{-3}$ & 2.08 &  & 7.05$\times10^{-3}$& 2.11 \\\hline
		$640^2$& 1.79$\times10^{-3}$ & 1.90 &  & 1.28$\times10^{-3}$ & 2.47 \\\hline
	\end{tabular}

\end{table}

\subsection{Non-linear equation in 3D: Translation}
We now consider the example in 3D. Interface $\Gamma_t$ is a moving sphere, which is represented as a zero level-set of function 
\[
\phi(x,y,z,t)= \sqrt{(x+0.5-t)^2+(y+0.5-t)^2+(z+0.5-t)^2}-1.
\]
For the exact solution 
\begin{equation*}
u(x,y,z,t)=\begin{cases}
1 & \text{in }\Omega^{-}\\
\frac{1}{\sqrt{(x+0.5-t)^2+(y+0.5-t)^2+(z+0.5-t)^2}} & \text{in }\Omega^{+}
\end{cases},
\end{equation*}
the velocity field is set as $V= \left(u,u,u \right)$. Since $u|_{\Gamma_t}=1$, movement of the interface $\Gamma_t$ agrees with $V$.
Numerical simulations are conducted up to final time $T=1$ with quantities $\rho^+=1, \rho^-=100$, and $\mu^+=1 ,\mu^-=10$. Convergence results and average iteration number of GMRES for each time step, denoted as $N_g$, are presented in table \ref{tab:ex6}. A second-order convergence of SL-GFM in 3D is verified; the iteration number does not dramatically increase.

\begin{table}[H]
	
	\centering
	\caption{Convergence rates for example 6\label{tab:ex6}}
	
	\begin{tabular}{|ccccccc|}
		\hline 
		Grid & $\parallel u(x,y,z,t)-u_{ijk}\parallel_{\infty}$ & order &  & $\parallel \phi(x,y,z,t)-\phi_{ijk}\parallel_{\infty}$ & order& $N_g$\tabularnewline
		\hline 
		$40^3$                 & 1.37$\times 10^{-2}$ &  -    &  & 1.47$\times 10^{-2}$ &  -  &16  \\ \hline
		$80^3$                 & 3.63$\times 10^{-3}$ & 1.92 &  & 3.70$\times 10^{-3}$ & 1.99 &23\\ \hline
		$160^3$                & 9.63$\times 10^{-4}$ & 1.91 &  & 9.60$\times 10^{-4}$ & 1.95 &32\\ \hline
	\end{tabular}

\end{table}

\section{Conclusion}
In this paper, a second-order accurate finite difference method for solving convection diffusion equations 
with jump conditions on a moving interface is presented. A bilinear interpolation using the ghost values near the interface is developed and adopted to the interpolation procedure of the semi-Lagrangian method. Coupled with second-order ghost fluid method \cite{cho2019second}, this produces a second-order convergence in $L^\infty$ norms.
	Furthermore, we have presented a second-order algorithm for an evolving interface when the velocity has jumps in its normal derivatives. Under the assumption that second-order time-discretization of two-phase incompressible flow is provided, we expect that the proposed method can be used to develop a second-order sharp capturing method for incompressible two-phase flows.
	\bibliographystyle{spmpsci} 

\bibliography{mybib}

\end{document}